\documentclass[leqno, myheadings, twoside]{amsart}
\usepackage{amsfonts}
\usepackage{amsmath, amsthm, amssymb, amscd, amsxtra,graphicx}
\usepackage{latexsym, amsfonts}
\usepackage{url}
\usepackage{texdraw}
\usepackage{epsfig}
\makeatletter \@addtoreset{equation}{section}

\makeatletter \renewcommand{\@biblabel}[1]{#1.}

\theoremstyle{remark}

\begin{document}
\title [The vibrational frequencies of the elastic body and its geometric quantities] {The vibrational frequencies of the elastic body and its geometric quantities}
\author{Genqian Liu}


\maketitle Department of Mathematics, Beijing Institute of
Technology,  Beijing 100081,  People's Republic of China.
 \ \    E-mail address:  liugqz@bit.edu.cn

\vskip 0.46 true cm

\vskip 0.45 true cm

\vskip 15 true cm

\begin{abstract}    For a bounded domain $\Omega\subset {\Bbb R}^n$ with smooth boundary, we explicitly calculate the first two coefficients of the asymptotic expansion of the trace of the strongly continuous semigroup associated with the Navier-Lam\'{e} operator on $\Omega$ as $t\to 0^+$. These coefficients (i.e., spectral invariants) provide precise information for the volume  of the elastic body $\Omega$ and the surface area of the boundary $\partial \Omega$ in terms of the spectrum of the  Navier-Lam\'{e} problem. As an application, we show that an $n$-dimensional ball is uniquely determined by its Navier-Lam\'{e} spectrum among all bounded elastic body with smooth boundary.    \end{abstract}

\vskip 1.39 true cm

\section{ Introduction}

\vskip 0.45 true cm

For the Navier-Lam\'{e} equations, one of important problems is to study the geometry of the elastic body from the
 vibrational frequencies of the elastic body, because this just reveals the true  behavior of the elastic body.
Let $\Omega\subset {\Bbb R}^n$ ($n\ge 2$) be
 a bounded domain with smooth boundary $\partial \Omega$.
 Let $P$ be the Navier-Lam\'{e} operator:
 \begin{eqnarray} \label {1-1} P\mathbf{u}:=-\tau \Delta \mathbf{u} -(\tau +\mu) \nabla (\nabla\cdot \mathbf{u}),   \quad \, \mathbf{u}=(u_1, u_2, \cdots, u_n),\end{eqnarray}
 where  $\nabla$  is the gradient  operator, $\Delta$ is the Laplacian, and $\tau$ and $\mu$ are Lam\'{e} parameters with $\tau>0$, $\tau+\mu>0$
 (see \cite{Ev}, \cite{Gur}, \cite{LaLi}, \cite{Sl},  \cite{Ar}).
    We denote by $P^-$ and $P^+$ the Navier-Lam\'{e} operators with the Dirichlet and Neumann boundary conditions, respectively.
   For the Derivation of the Navier-Lam\'{e} equation, its mechanical meaning and the explanation of the Dirichlet and Neumann boundary conditions, we refer the reader to \cite{Bant} or \cite{Sl}.

        According to theory of the elliptic equations, the Navier-Lam\'{e} operator $P$ can generate
   analytic semigroup $U^-(t)$ (respectively, $U^+(t)$)  with respect to the Dirichlet (respectively, Neumann) boundary condition in a suitable space of vector-valued functions (see  \cite{Mo1}, \cite{Mo2}, \cite{Mo3}, \cite{Pa}, \cite{St}). Furthermore,
   there exists a function-valued matrix ${\mathbf{K}}^- (t, x, y)$ (respectively, ${\mathbf{K}}^+(t,x,y)$), which is called the heat kernel or fundamental solution, such that
        \begin{eqnarray}  (U^-(t){\mathbf{v}}_0)(x)=\int_\Omega {\mathbf{K}}^-(t, x,y) {\mathbf{v}}_0(y)dy, \quad \,
        {\mathbf{v}}_0\in  [H^1_0(\Omega)]^n\end{eqnarray}
        (respectively,  \begin{eqnarray}  (U^+(t){\mathbf{v}}_0)(x)=\int_\Omega {\mathbf{K}}^+(t, x,y) {\mathbf{v}}_0(y)dy, \quad \,
        {\mathbf{v}}_0\in  [H^1(\Omega)]^n).\end{eqnarray}
        Thus, ${\mathbf{v}}^-(t, x):= (U^-(t) \mathbf{v}_0)(x)$ and ${\mathbf{v}}^+(t, x):= (U^+(t) \mathbf{v}_0)(x)$ respectively satisfy the initial-boundary problems for the elastodynamic evolution equations
          \begin{eqnarray} \label{1-2}  \left\{ \begin{array}{ll}  {\mathbf{v}}^-_t -
          \tau \Delta {\mathbf{v}}^- -(\tau +\mu) \nabla (\nabla\cdot {\mathbf{v}}^-)=0
           \;\; &\mbox{in}\;\;  (0,+\infty)\times \Omega,\\
            {\mathbf{v}}^-=0 \;\;& \mbox{on}\;\; (0,+\infty)\times \partial \Omega, \\
     {\mathbf{v}}^-(0,x)= \mathbf{v}_0 \;\; & \mbox{on}\;\; \{0\}\times \Omega \end{array}  \right.\end{eqnarray}
       and   \begin{eqnarray} \label{1-3}  \left\{ \begin{array}{ll}  {\mathbf{v}}^+_t -
          \tau \Delta {\mathbf{v}}^+ -(\tau +\mu) \nabla (\nabla\cdot {\mathbf{v}}^+)=0
           \;\; &\mbox{in}\;\;  (0,+\infty)\times \Omega,\\
            \frac{\partial {\mathbf{v}}^+}{\partial \nu}=0 \;\;& \mbox{on}\;\; (0,+\infty)\times \partial \Omega, \\
     {\mathbf{v}}^+(0,x)= \mathbf{v}_0 \;\; & \mbox{on}\;\; \{0\}\times \Omega. \end{array}  \right.\end{eqnarray}

     On the other hand, since the Navier-Lam\'{e} operator $P^-$ (respectively, $P^+$) is a unbounded self-adjoint positive operator
     in $[H^1_0(\Omega)]^n$ (respectively, $(H^1(\Omega)$) with discrete spectrum $0< \lambda_1^- < \lambda_2^- \le \cdots \le \lambda_k^- \le \cdots \to +\infty$ (respectively, $0=\lambda_1^+ < \lambda_2^+ \le \cdots \le \lambda_k^+ \le \cdots \to +\infty$),
 one has \begin{eqnarray} \label{1-4} P^\mp {\mathbf{u}}_k^\mp =\lambda_k^\mp {\mathbf{u}}_k^\mp,\end{eqnarray} where ${\mathbf{u}}_k^-\in  [H^1_0(\Omega)]^n$ (or ${\mathbf{u}}_k^+\in  [H^1(\Omega)]^n$) are the corresponding orthogonal eigenvectors.
(\ref{1-4}) can be rewritten as
  \begin{eqnarray} \label{1-5} \left\{ \begin{array}{ll}   -\tau \Delta {\mathbf{u}}_k^- - (\tau+\mu) \nabla (\nabla\cdot {\mathbf{u}}_k^-)  = \lambda^-_k {\mathbf{u}}_k^-  \;\; &\mbox{in}\;\;  \Omega,\\
        {\mathbf{u}}_k^-=0 \;\;& \mbox{on}\;\;  \partial \Omega
    \end{array}  \right.\end{eqnarray}
    and
    \begin{eqnarray} \label{1-6} \left\{ \begin{array}{ll}   -\tau \Delta {\mathbf{u}}_k^+ - (\tau+\mu) \nabla (\nabla\cdot {\mathbf{u}}_k^+)  = \lambda^+_k {\mathbf{u}}_k^+  \;\; &\mbox{in}\;\;  \Omega,\\
       \frac{\partial  {\mathbf{u}}_k^+}{\partial \nu}=0 \;\;& \mbox{on}\;\;  \partial \Omega.
    \end{array}  \right.\end{eqnarray}
     Clearly, the eigenvalue problem (\ref{1-5}) and (\ref{1-6}) can be immediately obtained by considering the solution of the form ${\mathbf{v}}^\mp(t, x)=T(t) {\mathbf{u}}^\mp(x)$ in Navier-Lam\'{e} evolution equations (\ref{1-2}) and (\ref{1-3}). The Navier-Lam\'{e} eigenvalues are physical quantities because
  they just are the vibrational frequencies of an elastic body in the two or three dimensions.

 An interesting question, which is similar to the famous Kac question for the Dirichlet-Laplacian
  (see \cite{Kac}, \cite{Lo} or \cite{We1} ), is: ``can one hear the shape of an elastic body by hearing the vibrational frequencies (or pitches) of elastic body ?''
More precisely, we have elastic body of different shapes. You hit them, and then you listen to the
 frequencies of elastic wave. Can you tell the shape (or the geometric quantities) of the elastic body?

In this paper, some surprising and interesting results are obtained by considering the Navier-Lam\'{e}
operator semigroup $U^\mp(t) = e^{-tP^\mp}$ and by using some new methods of pseudodifferential
operators. The following theorem is the main result of this paper:

\vskip 0.25 true cm

  \noindent{\bf Theorem 1.1.} \ {\it Let  $\Omega\subset {\Bbb R}^n$ ($n\ge 2$) be
 a bounded domain with smooth boundary $\partial \Omega$, and let $0< \lambda_1^-< \lambda_2^- \le \cdots \le \lambda_k^- \le \cdots$
 (respectively, $0= \lambda_1^+ < \lambda_2^+ \le \lambda_3^+ \le \cdots \le \lambda_k^+ \le \cdots $) be the eigenvalues
 of the Navier-Lam\'{e} operator $P^-$ (respectively, $P^+$) with respect to the Dirichlet (respectively, Neumann) boundary condition.  Then
 \begin{eqnarray} \label{1-7} &&  \sum_{k=1}^\infty e^{-\lambda_k^\mp t}  = \mbox{Tr}(e^{-tP^\mp})=\bigg[ \frac{(n-1)}{(4\pi \tau t)^{n/2}}
 + \frac{1}{(4\pi (2\tau+\mu) t)^{n/2}}\bigg] |\Omega| \\
&&  \quad \, \mp \frac{1}{4} \bigg[  \frac{(n-1)}{(4\pi \tau t)^{(n-1)/2}}
 +  \frac{1}{(4\pi (2\tau+\mu) t)^{(n-1)/2}}\bigg]|\partial \Omega| +O(t^{1-\frac{n}{2}})\quad\;\; \mbox{as}\;\; t\to 0^+.\nonumber\end{eqnarray}
 Here $|\Omega|$ denotes the $n$-dimensional volume of $\Omega$, and $|\partial \Omega|$ denotes the $(n-1)$-dimensional volume of  $\partial \Omega$.}

 \vskip 0.25 true cm

  Our result shows that not only the volume $|\Omega|$  but also the surface area $|\partial \Omega|$ can be known if  we know all Navier-Lam\'{e} eigenvalues with respect to the Dirichlet (respectively, Neumann) boundary condition. Roughly speaking,  one can ``hear'' the volumes of the domain and the surface area of its boundary $\partial \Omega$ by ``hearing'' all the pitches of the vibration of an elastic body.

The key ideas of this paper are as follows. If ${\mathbf{u}}_k^\mp$ is the normalized eigenvector of Navier-Lam\'{e} problem with eigenvalue $\lambda_k^\mp$, the Navier-Lam\'{e} heat kernel ${\mathbf{K}}^\mp(t, x, y)$ is  given by \begin{eqnarray} \label{1-0a-1} {\mathbf{K}}^\mp(t,x,y) =\sum_{k=1}^\infty e^{-t \lambda_k^\mp} {\mathbf{u}}_k^\mp(x)\otimes {\mathbf{u}}_k^\mp(y).\end{eqnarray}
Thus the integral of the trace of ${\mathbf{K}}^\mp(t,x,y)$ is actually a spectral invariants: by  (\ref{1-0a-1}), we can compute
\begin{eqnarray} \label{1-0a-2} Tr\bigg(\int_\Omega {\mathbf{K}}^\mp(t,x,x) dx\bigg)=\sum_{k=1}^\infty e^{-t \lambda_k^\mp}.\end{eqnarray}
To further analyze the geometric content of the spectrum, we calculate the same
trace by an entirely different way: we constructs the heat kernel
 from $P^\pm$ by the Cauchy integral formula:
 \begin{eqnarray*} e^{-tP^\mp} =\frac{1}{2\pi i} \int_{\mathcal{C}} e^{-t\lambda} (\lambda I- P^\mp)^{-1} d\lambda,\end{eqnarray*}
where $\mathcal{C}$ is a suitable curve in the complex plane in the positive direction around the spectrum of $P^\mp$.
  By calculating the full symbols of pseudodifferential operator, and then applying technique of Mckean-Singer, we explicitly computes the integral of the
trace of $e^{-tP^\mp}$. We can show  that the integral of the trace has an asymptotic expansion 
\begin{eqnarray} Tr\bigg(\int_{\Omega} {\mathbf{K}}^\mp(t,x,x) dx\bigg) \sim a_0t^{-n/2}+ a_1^\mp t^{-(n-1)/2} +\cdots\quad \quad \mbox{as}\;\; t\to 0^+,\end{eqnarray}
   where $a_0=\big[ \frac{(n-1)}{(4\pi \tau t)^{n/2}}
 + \frac{1}{(4\pi (2\tau+\mu) t)^{n/2}}\big] |\Omega|$, $ \;a_1^\mp=\mp \frac{1}{4}\big[  \frac{(n-1)}{(4\pi \tau t)^{(n-1)/2}}
 +  \frac{1}{(4\pi (2\tau+\mu) t)^{(n-1)/2}}\big]|\partial \Omega|$. More exactly,
   we can apply the Seeley's calculus in the interior of $\Omega$ (see \cite{See} or \cite{Gr}) to the  symbols of the Navier-Lam\'{e} operator to get the  coefficient $a_0$.  However, the Seeley's method can't be used to deal with the boundary case for the Navier-Lam\'{e} operator (cf. \cite{See} or \cite{Gre}). This has become a stumbling block in the expansion of the heat trace for the corresponding elastic operator. To overcome this problem and to obtain the second coefficient $a_1$, we will approximate the heat kernel near the boundary locally by the ``method of images.''  Locally, the boundary looks like the superplane $x_n = 0$ in the ${\Bbb R}^n$; letting $x\to \overset{*}{x}$ be
the reflection $(x_1, \cdots, x_{n-1}, x_n) \to (x_1, \cdots, x_{n-1},-x_n)$, the kernel
$ {\mathbf{K}}^\mp(t, x, y) = {\mathbf{K}}(t, x, y) \mp {\mathbf{K}}(t, x, \overset{*}{y})$ (or its normal derivative) vanishes on $x_n=0$. By further estimating the traces of these two heat kernels, we finally obtain coefficient $a_1^\mp$.

\vskip 0.10 true cm

\vskip 0.15 true cm

As an application of theorem 1.1, we can prove the following spectral rigidity result:
\vskip 0.25 true cm

  \noindent{\bf Corollary 1.2.} \ {\it Let $\Omega \subset {\Bbb R}^3$ be a bounded domain with smooth boundary $\partial \Omega$.
   Suppose that the Navier-Lam\'{e} spectrum with respect to the Dirichlet (respectively, Neumann) boundary condition, is equal to that of $B_r$, a ball of radius $r$. Then $\Omega=B_r$. }

\vskip 0.23 true cm

  Corollary 1.2 also shows that a ball is uniquely determined by its Navier-Lam\'{e} spectrum among all Euclidean bounded domains with smooth boundary.

\vskip 0.25 true cm

\vskip 1.49 true cm

\section{Some notations and a Lemmas}

\vskip 0.45 true cm

If $W$ is an open subset  of ${\Bbb R}^n$, we denote  by $S^m_{1,0}=S^m_{1,0} (W,
{\Bbb R}^n)$ the set of all $p\in C^\infty (W, {\Bbb R}^n)$ such
that for every compact set $O\subset W$ we have
 \begin{eqnarray} \label {-2.3} |D^\beta_x D^\alpha_\xi p(x,\xi)|\le C_{O,\alpha,
 \beta}(1+|\xi|)^{m-|\alpha|}, \quad \; x\in O,\,\, \xi\in {\Bbb R}^n\end{eqnarray}
 for all $\alpha, \beta\in {\Bbb N}^n_+$.
 The  elements  of $S^m_{1,0}$  are  called  symbols (or full symbols) of order $m$.
 It is clear that $S^m_{1,0}$ is a
Fr\'{e}chet space with semi-norms given by the smallest constants
which can be used in (\ref{-2.3}) (i.e.,
\begin{eqnarray*} \|p\|_{O,\alpha, \beta}=
 \,\sup_{x\in O}\big|\left(D_x^\beta D_\xi^\alpha
 p(x, \xi)\right)(1+|\xi|)^{|\alpha|-m}\big|).\end{eqnarray*}
 Let $p(x, \xi)\in S^m_{1,0}$. A pseudo-differential operator
in an open set $W\subset {\Bbb R}^n$ is essentially defined by a
Fourier integral  operator (cf. \cite{KN}, \cite{Ho4}, \cite{Ta2}, \cite{Gr}):
\begin{eqnarray} \label{-2.1}  P(x,D) u(x) = \frac{1}{(2\pi)^{n}} \int_{{\Bbb R}^n} p(x,\xi)
 e^{i \langle x,\xi\rangle} \hat{u} (\xi)d\xi,\end{eqnarray}
and denoted by $OPS^m$.
Here $u\in C_0^\infty (W)$ and $\hat{u} (\xi)=
 \int_{{\Bbb R}^n} e^{-i\langle y, \xi\rangle} u(y)dy$ is the Fourier
transform of $u$.
 If there are smooth $p_{m-j} (x, \xi)$, homogeneous in $\xi$ of degree $m-j$ for $|\xi|\ge 1$, that is,
  $p_{m-j} (x, r\xi) =r^{m-j} p_{m-j} (x, \xi)$ for $r, |\xi|\ge 1$, and if
 \begin{eqnarray} \label{-2.2} p(x, \xi) \sim \sum_{j\ge 0} p_{m-j} (x, \xi)\end{eqnarray}
 in the sense that
 \begin{eqnarray}  p(x, \xi)-\sum_{j=0}^l p_{m-j} (x, \xi) \in S^{m-l-1}_{1, 0}, \end{eqnarray}
 for all $l$, then we say $p(x,\xi) \in S_{cl}^m$, or just $p(x, \xi)\in S^m$. We call
 $p_m(x, \xi)$ the principal symbols of $P(x, D)$.

\vskip 0.16 true cm

 An operator $P$ is said to be an elliptic pseudodifferential operator of order $m$ if
  for every compact $O\subset \Omega$ there exists a  positive constant $c=c(O)$ such that \begin{eqnarray*} |p(x, \xi)\ge c|\xi|^m, \,x\in O,\, |\xi|\ge 1\end{eqnarray*} for any compact set $O\subset \Omega$.
  If $P$ is a non-negative elliptic pseudodifferential operator of order $m$, then the spectrum  of $P$ lies in a right half-plane and has a finite lower bound $\tau(P) =
\inf\{\mbox{Re}\, \lambda\big| \lambda\in \sigma(P)\}$. We can modify $p_m(x, \xi)$ for small $\xi$ such that
$p_m(x,\xi)$  has a positive
lower bound throughout and lies in $\{\lambda=re^{i\theta} \big| r>0, |\theta|\le \theta_0\}$, where $\theta_0\in (0,\frac{\pi}{2})$.
According to \cite{Gr}, the resolvent $(P-\lambda)^{-1}$ exists and is holomorphic in $\lambda$ on a neighborhood of a set
\begin{eqnarray*}  W_{r_0,\epsilon} =\{\lambda\in {\Bbb C} \big| |\lambda|\ge r_0, \mbox{arg}\, \lambda \in [\theta_0+\epsilon, 2\pi -\theta_0-\epsilon], \,\mbox{Re}\, \lambda \le \tau (P)-\epsilon\}\end{eqnarray*}
(with $\epsilon>0$). There exists a parametrix $Q'_\lambda$ on a neighborhood of a possibly larger set
(with $\delta>0,\epsilon>0$)
\begin{eqnarray*} V_{\delta,\epsilon} =\{ \lambda \in {\Bbb C} \big| |\lambda| \ge \delta \;\;\mbox{or arg}\, \lambda \in [\theta_0+\epsilon, 2\pi-\theta_0-\epsilon]\}\end{eqnarray*}
such that this parametrix coincides with $(P-\lambda)^{-1}$ on the intersection. Its symbol $q(x,\xi, \lambda)$
in local coordinates is holomorphic in $\lambda$ there and has the form (cf. Section 3.3 of \cite{Gr})
\begin{eqnarray} \label {-2.6} q(x, \xi,\lambda) \sim \sum_{l\ge 0} q_{-m-l} (x, \xi, \lambda)\end{eqnarray}
where \begin{eqnarray} \label{-2.7} & q_{-m}= (p_m(x, \xi) -\lambda)^{-1}, \quad \;  q_{-m-1} = b_{1,1}(x, \xi) q^2_{-m},\\
& \, \cdots, \,  q_{-m-l} = \sum_{k=1}^{2l} b_{l,k} (x, \xi) q^{k+1}_{-m}, \cdots, \nonumber\end{eqnarray}
with symbols $b_{l,k}$ independent of $\lambda$  and homogeneous of degree $mk-l$ in $\xi$ for $|\xi|\ge1$.
 The semigroup $e^{-tP}$ can be defined from $P$ by the Cauchy integral formula (see p.$\,$4 of \cite{GG}):
 \begin{eqnarray*} e^{-tP} =\frac{1}{2\pi i} \int_{\mathcal{C}} e^{-t\lambda} (\lambda I -P)^{-1} d\lambda,\end{eqnarray*}
where $\mathcal{C}$ is a suitable curve in the complex plane in the positive direction around the spectrum of $P$.

\vskip 0.15 true cm

Clearly, the Navier-Lam\'{e} operator $ P\mathbf{u}=-\tau \Delta \mathbf{u} -(\tau +\mu) \nabla (\nabla\cdot \mathbf{u}),  \quad  \mathbf{u}= (u_1, u_2, \cdots, u_n)$,
can be rewritten as $P\mathbf{u}(x) =\frac{1}{(2\pi)^n}\int_{{\Bbb R}^n} ({\mathbf{A}} (x,\xi))e^{i\langle x,\xi\rangle} {\hat{\mathbf{u}}}^T(\xi)\,d\xi$,
where   \begin{gather} \label{3-3-00}
   \mathbf{A}(x, \xi) =\begin{pmatrix} \tau |\xi|^2 +(\tau+\mu)\xi_1^2   & (\tau+\mu) \xi_1\xi_2 & \cdots \cdots  & (\tau+\mu)\xi_1\xi_n\\
    (\tau+\mu)\xi_2\xi_1 & \tau |\xi|^2+ (\tau+\mu) \xi_2^2  & \cdots\cdots & (\tau+\mu) \xi_2 \xi_n\\
    \cdots \cdots\cdots \cdots  \cdots&\cdots \cdots\cdots \cdots \cdots&\cdots \cdots&\cdots\cdots \cdots\cdots \\
    (\tau+\mu)\xi_n\xi_1 & (\tau+\mu) \xi_n\xi_2& \cdots\cdots  & \tau|\xi|^2 +(\tau+\mu)\xi_n^2 \end{pmatrix} \end{gather}
and $\hat{\mathbf{u}}^T= ({\hat{u}}_1, {\hat{u}}_2,\cdots, {\hat{u}}_n)^T$.

\vskip 0.18 true cm

 We will calculate the reverse of matrix $\lambda I - \mathbf{A}$ and further give its trace.

\vskip 0.18 true cm

  \noindent{\bf Lemma 2.1.} \ {\it  For any integers $n\ge 1$, we have
   \begin{eqnarray} \mbox{Tr}\big((\lambda I- \mathbf{A})^{-1}\big)&=&(\lambda -\tau |\xi|^2)^{-1}
  (\lambda-(2\tau+\mu)|\xi|^2)^{-1}  \\
   &&   \times \big[ n\lambda-
  \big((2n-1)\tau +(n-1)\mu\big)|\xi|^2\big].\nonumber\end{eqnarray}}

 \vskip 0.28 true cm

  \noindent  {\it Proof.} \  (i) \ When $n=1$, it is clear that
  \begin{eqnarray*} (\lambda \mathbf{I}- \mathbf{A})^{-1} = (\lambda \mathbf{I}- (2\tau+\mu)|\xi|^2)^{-1}.\end{eqnarray*}

\vskip 0.15 true cm

(ii) \  When $n=2$, we have
  \begin{eqnarray*}   (\lambda I- \mathbf{A})^{-1} &=& (\lambda -\tau |\xi|^2)^{-1} \big(\lambda -(2\tau+\mu)|\xi|^2\big)^{-1} \\
  && \; \times \begin{pmatrix} \lambda -\tau |\xi|^2 -(\tau+\mu) \xi_2^2    & (\tau+\mu) \xi_1\xi_2 \\
     (\tau+\mu)\xi_2\xi_1 &   \lambda -\tau|\xi|^2 -(\tau+\mu) \xi_1^2\end{pmatrix},
 \end{eqnarray*} which implies
  \begin{eqnarray*} \mbox{Tr} \big((\lambda I- \mathbf{A})^{-1} \big)= (\lambda -\tau |\xi|^2)^{-1}(\lambda -(2\tau+\mu)|\xi|^2)^{-1} (2\lambda -(3\tau +\mu)|\xi|^2).\end{eqnarray*}

\vskip 0.15 true cm

  (iii) \ \  When $n=3$, we get  \begin{eqnarray*} (\lambda I- \mathbf{A})^{-1} = (\lambda -\tau |\xi|^2)^{-2}
   (\lambda -(2\tau +\mu) |\xi|^2\big)^{-1}\cdot (\lambda -\tau|\xi|^2)  \qquad \qquad \qquad \quad\qquad\qquad \qquad \qquad\\
\times \begin{pmatrix}  \lambda -\tau |\xi|^2 -(\tau+\mu) (\xi_2^2 +\xi_3^2) &
    (\tau+\mu)\xi_1\xi_2 &   (\tau+\mu)\xi_1\xi_3\\
  (\tau+\mu) \xi_2\xi_1 &  \lambda -\tau |\xi|^2 -(\tau+\mu) (\xi_1^2 +\xi_3^2) &
    (\tau+\mu) \xi_2\xi_3\\
    (\tau+\mu) \xi_3\xi_1 & (\tau+\mu) \xi_3\xi_2 & \lambda -\tau |\xi|^2 -(\tau+\mu) (\xi_1^2 +\xi_2^2)\end{pmatrix}\end{eqnarray*}
 so that
  \begin{eqnarray*} \mbox{Tr} (\lambda \mathbf{I}- \mathbf{A})^{-1} = (\lambda -\tau |\xi|^2)^{-1}(\lambda -(2\tau+\mu)|\xi|^2)^{-1} (3\lambda -(5\tau +2\mu)|\xi|^2).\end{eqnarray*}

\vskip 0.15 true cm

(iv) \ \  For $n=4$, we find that
 \begin{eqnarray*} (\lambda I- \mathbf{A})^{-1} = (\lambda -\tau |\xi|^2)^{-3}
   (\lambda -(2\tau +\mu) |\xi|^2\big)^{-1} \begin{pmatrix} b_{11}^* & b_{12}^* & b_{13}^* & b_{14}^*\\
 b_{21}^* & b_{22}^* & b_{23}^* & b_{24}^*\\
 b_{31}^* & b_{32}^* & b_{33}^* & b_{34}^*\\
 b_{41}^* & b_{42}^* & b_{43}^* & b_{44}^* \end{pmatrix},\end{eqnarray*}
 where  $b_{ij}^*$ is the (i,j)-cofactor of the matrix $(\lambda- A)$.  It is easy to check that \begin{eqnarray*} b_{11}^*&=& \big(\lambda -\tau |\xi|^2 -(\tau+\mu) \xi_2^2\big)\big(\lambda -\tau |\xi|^2 -(\tau+\mu) \xi_3^2\big)\big(\lambda -\tau |\xi|^2 -(\tau+\mu) \xi_4^2\big)\\
  && -2 (\tau+\mu)^3 \xi_2^2\xi_3^2 \xi_4^2 -\big(\lambda -\tau |\xi|^2 -(\tau+\mu) \xi_2^2\big)(\tau+\mu)^2 \xi_3^2 \xi_4^2\\
   && -\big(\lambda -\tau |\xi|^2 -(\tau+\mu) \xi_3^2\big)(\tau+\mu)^2 \xi_2^2 \xi_4^2
   - \big(\lambda -\tau |\xi|^2 -(\tau+\mu) \xi_4^2\big)(\tau+\mu)^2 \xi_2^2 \xi_3^2,\end{eqnarray*}
\begin{eqnarray*} b_{22}^*&=& \big(\lambda -\tau |\xi|^2 -(\tau+\mu) \xi_1^2\big)\big(\lambda -\tau |\xi|^2 -(\tau+\mu) \xi_3^2\big)\big(\lambda -\tau |\xi|^2 -(\tau+\mu) \xi_4^2\big)\\
  && -2 (\tau+\mu)^3 \xi_1^2\xi_3^2 \xi_4^2  -\big(\lambda -\tau |\xi|^2 -(\tau+\mu) \xi_1^2\big)(\tau+\mu)^2 \xi_3^2 \xi_4^2\\
  && -\big(\lambda -\tau |\xi|^2 -(\tau+\mu) \xi_3^2\big)(\tau+\mu)^2 \xi_1^2 \xi_4^2
      - \big(\lambda -\tau |\xi|^2 -(\tau+\mu) \xi_4^2\big)(\tau+\mu)^2 \xi_1^2 \xi_3^2,\end{eqnarray*}
\begin{eqnarray*} b_{33}^*&=& \big(\lambda -\tau |\xi|^2 -(\tau+\mu) \xi_1^2\big)\big(\lambda -\tau |\xi|^2 -(\tau+\mu) \xi_2^2\big)\big(\lambda -\tau |\xi|^2 -(\tau+\mu) \xi_4^2\big)\\
  && -2 (\tau+\mu)^3 \xi_1^2\xi_2^2 \xi_4^2 -\big(\lambda -\tau |\xi|^2 -(\tau+\mu) \xi_1^2\big)(\tau+\mu)^2 \xi_2^2 \xi_4^2 \\
  && -\big(\lambda -\tau |\xi|^2 -(\tau+\mu) \xi_2^2\big)(\tau+\mu)^2 \xi_1^2 \xi_4^2- \big(\lambda -\tau |\xi|^2 -(\tau+\mu) \xi_4^2\big)(\tau+\mu)^2 \xi_1^2 \xi_2^2,\end{eqnarray*}
\begin{eqnarray*} b_{44}^*&=& \big(\lambda -\tau |\xi|^2 -(\tau+\mu) \xi_1^2\big)\big(\lambda -\tau |\xi|^2 -(\tau+\mu) \xi_2^2\big)\big(\lambda -\tau |\xi|^2 -(\tau+\mu) \xi_3^2\big)\\
  && -2 (\tau+\mu)^3 \xi_1^2\xi_2^2 \xi_3^2 -\big(\lambda -\tau |\xi|^2 -(\tau+\mu) \xi_1^2\big)(\tau+\mu)^2 \xi_2^2 \xi_3^2 \\
  && -\big(\lambda -\tau |\xi|^2 -(\tau+\mu) \xi_2^2\big)(\tau+\mu)^2 \xi_1^2 \xi_3^2- \big(\lambda -\tau |\xi|^2 -(\tau+\mu) \xi_3^2\big)(\tau+\mu)^2 \xi_1^2 \xi_2^2.\end{eqnarray*}
 A direct calculation shows that
 \begin{eqnarray*} b_{11}^* + b_{22}^* + b_{33}^* +b_{44}^* = (\lambda - \tau|\xi|^2 )^2 (4\lambda -(7\tau +3\mu)|\xi|^2).\end{eqnarray*}
 Thus, we have that
 \begin{eqnarray*} \mbox{Tr}  \big( (\lambda \mathbf{I}-\mathbf{A})^{-1}\big)&=& \big(\lambda-\tau |\xi|^2)^{-3} (\lambda -(2\tau+\mu)|\xi|^2)^{-1} (\lambda -\tau|\xi|^2 )^2 \big(4\lambda -(7\tau +3\mu)|\xi|^2\big)\\
 &=& (\lambda- \tau |\xi|^2)^{-1} \big(\lambda -(\tau+\mu)|\xi|^2\big)^{-1} \big(4\lambda -(7\tau +3\mu)|\xi|^2\big).\end{eqnarray*}

\vskip 0.15 true cm

 (v) \ \  For general integer $n\ge 1$, since  \begin{eqnarray*} (\lambda \mathbf{I} - \mathbf{A})^{-1}= |\lambda \mathbf{I} - \mathbf{A}|^{-1} \begin{pmatrix} c_{11}^*  & c_{12}^*  & \cdots & c_{1n}^*\\
 c_{21}^* & c_{22}^* & \cdots & c_{2n}^*\\
 \cdots  & \cdots  & \cdots   &  \cdots \\
 c_{n1}^* & c_{n2}^* & \cdots & c_{nn}^* \end{pmatrix},\end{eqnarray*}
 where $c_{ij}^*$ is the (i,j)-cofactor of the matrix $(\lambda {\mathbf{I}}- \mathbf{A})$.
 It follows from the technique of (iv) that for any integers $n\ge 1$,
  \begin{eqnarray*} |\lambda \mathbf{I}- \mathbf{A}|=(\lambda -\tau |\xi|^2)^{(n-1)} \big(\lambda -(2\tau+\mu)|\xi|^2\big)\end{eqnarray*}
  and \begin{eqnarray*} c_{ii}^* &=& \prod_{\underset {j\ne i}{1\le j\le n}} \big( \lambda -\tau|\xi|^2 -(\tau+\mu)\xi_j^2\big)
   - \sum_{k=0}^{n-3} (n-k-2) \big(\lambda -\tau|\xi|^2-(\tau+\mu) \xi^2_{l_1}\big )\\ && \times \cdots \times\big(\lambda -\tau|\xi|^2-(\tau+\mu) \xi^2_{l_k}\big )  (\tau+\mu)^{n-k-1} \xi^2_{m_1} \cdots \xi^2_{m_{n-k-1}}, \quad \, i= 1,2, \cdots, n,\end{eqnarray*}
   where \begin{eqnarray*} & l_1, \cdots, l_k\in \{ 1,2, \cdots, n\}\setminus \{i\}, \quad \; l_1< \cdots < l_k, \;\;\mbox{and}\\
     & m_1, \cdots, m_{n-k-1}\in \{ 1,2, \cdots, n\}\setminus \{i, l_1, \cdots, l_k\}, \,\quad m_1< \cdots < m_{n-k-1}.\end{eqnarray*}
   It follows that
   \begin{eqnarray*} \sum_{i=1}^n c_{ii}^* = (\lambda- \tau|\xi|^2)^{n-2} \big[n\lambda -\big((2n-1)\tau +(n-1)\mu\big)|\xi|^2\big].\end{eqnarray*}
      Hence
     \begin{eqnarray*}  \mbox{Tr}\big((\lambda \mathbf{I}- \mathbf{A})^{-1}\big)&=& (\lambda -\tau |\xi|^2)^{-(n-1)}
  (\lambda-(2\tau+\mu)|\xi|^2)^{-1} (\lambda -\tau |\xi|^2)^{(n-2)} \\
 &&\; \times \big[ n\lambda-
  \big((2n-1)\tau +(n-1)\mu\big)|\xi|^2\big]\nonumber\\
  &=&  (\lambda -\tau |\xi|^2)^{-1}
  (\lambda-(2\tau+\mu)|\xi|^2)^{-1}\big[ n\lambda-
  \big((2n-1)\tau +(n-1)\mu\big)|\xi|^2\big].\nonumber\end{eqnarray*} $\;\; \square$

\vskip 0.18 true cm

 Now, we construct a pseudodifferential operator $B$ to approximate the resolvant
$(\lambda I -P)^{-1}$ as follows: let \begin{eqnarray} \label{3--1} {\mathbf{b}}(x, \xi, \lambda) \sim {\mathbf{b}}_0(x, \xi, \lambda) +\cdots + {\mathbf{b}}_m(x, \xi, \lambda) +\cdots\end{eqnarray} be the expansion of the full symbol of $B$.
 Suppose that the complex parameter $\lambda$ have homogeneity $2$ (cf. \cite{See} or \cite{Gil}). Let the ${\mathbf{b}}_m$ be homogeneous
of order $-2-m$ in the variables $(\xi, \lambda)$. This infinite sum defines $\mathbf{b}$ asymptotically. Our purpose is to define
 $\mathbf{b}$ so that
    \begin{eqnarray}\label{3/1}  \sigma(B(\lambda I- P))\sim I,\end{eqnarray}
 where $\sigma(T)$ denotes the full symbol of pseudodifferential operator $T$.

  \vskip 0.10 true cm

We have the following:

 \vskip 0.18 true cm

  \noindent{\bf Lemma 2.2.} \ {\it  Let the full symbol $\mathbf{b}$ of $B$ have the form (\ref{3--1}) and let $B$ satisfy (\ref{3/1}) (i.e.,
 \begin{eqnarray}\label {3/2}  \sum_{\alpha\ge 0} (d_\xi^\alpha \mathbf{b})\cdot (D_x^\alpha (\sigma (\lambda I- P)))/\alpha !\,\sim I\,\,). \end{eqnarray}
  Then  ${\mathbf{b}}_0(x, \xi, \lambda)= (\lambda \mathbf{I}- \mathbf{A})^{-1}$ and ${\mathbf{b}}_m(x, \xi, \lambda)=0$ for all $m\ge 1$}.

  \vskip 0.18 true cm

    \noindent  {\it Proof.} \ Denote ${\mathbf{a}}_2 =(\lambda I- \mathbf{A})$, ${\mathbf{a}}_1=0$ and ${\mathbf{a}}_0=0$. We decompose  sum (\ref{3/1}) into orders of homogeneity (see \cite{Ho3}, \cite{Ho4} or p.$\,$13 of \cite{Ta2}):
 \begin{eqnarray*}  \sigma (B(\lambda I- P))\sim \sum_{m=0}^\infty \bigg( \sum_{m=j+|\alpha|+2-k} (d_\xi^\alpha {\mathbf{b}}_j )\cdot (D_x^\alpha {\mathbf{a}}_k)/\alpha !\bigg) .\end{eqnarray*} The sum is over terms which are homogeneous of order $-m$.
    Thus (\ref{3/1}) leads to the following equations
    \begin{eqnarray*} I &=& \sum_{0=j+|\alpha|+2-k} (d_\xi^\alpha {\mathbf{b}}_j) (D_x^\alpha {\mathbf{a}}_k)/\alpha ! ={\mathbf{b}}_0 {\mathbf{a}}_2,\\
    0&=& \sum_{m=j+|\alpha|+2-k} (d^\alpha_\xi {\mathbf{b}}_j)(D_x^\alpha {\mathbf{a}}_k)/\alpha !\\
    &=& {\mathbf{b}}_m {\mathbf{a}}_2 +  \sum_{\underset {j<m} {m=j+|\alpha|+2-k}} (d^\alpha_\xi {\mathbf{b}}_j)(D_x^\alpha {\mathbf{a}}_k)/\alpha !,\quad\quad m\ge 1.\end{eqnarray*}
     These equations determine the ${\mathbf{b}}_m$ inductively. In other words,
     we have \begin{eqnarray*}&& {\mathbf{b}}_0 ={\mathbf{a}}_2^{-1}, \,\quad   {\mathbf{b}}_m = -{\mathbf{b}}_0 \big( \sum_{j<m} (d^\alpha_\xi 
     {\mathbf{b}}_j )(D^\alpha_x {\mathbf{a}}_k)/\alpha !\big),\\
       &&   \qquad \qquad \qquad \qquad \qquad \qquad \mbox{for}\;\; m= j+|\alpha|+2-k.\end{eqnarray*}
    Since ${\mathbf{a}}_1={\mathbf{a}}_0=0$ and since ${\mathbf{a}}_2$ is independent of $x$, it can be immediately see that ${\mathbf{b}}_m(x, \xi,\lambda)=0$ for all $m\ge 1$. $\quad \quad \square$

\vskip 1.29 true cm

\section{Asymptotic expansion}

\vskip 0.45 true cm

 \noindent  {\it Proof of Theorem 1.1.} \  We calculate the asymptotic expansion of the trace of semigroup $e^{-tP}$ as $t\to 0^+$.  Note that the interior asymptotics are independent of the boundary condition; however, the
boundary asymptotics depend on the Dirichlet (or Neumann) boundary conditions. It follows from Lemmas 2.1 and 2.2 that the full symbols of the operator $(\lambda I- P)^{-1}$ is $(\lambda \mathbf{I}- \mathbf{A})^{-1}$. In view of
  \begin{eqnarray*} e^{-tP} =\frac{1}{2\pi i} \int_{\mathcal{C}} e^{-t\lambda} (\lambda I-P)^{-1} d\lambda,\end{eqnarray*}
where $\mathcal{C}$ is a suitable curve in the complex plane in the positive direction around the spectrum of $P$,
 we find that \begin{eqnarray}\quad  {\mathbf{K}} (t,x,y) = e^{-tP}\delta(x-y) = \frac{1}{(2\pi)^n} \int_{{\Bbb R}^n} e^{i\langle x-y, \xi\rangle} \bigg[\frac{1}{2\pi i} \int_{\mathcal{C}} e^{-t\lambda}  (\lambda \mathbf{I} -\mathbf{A})^{-1} d\lambda\bigg] d\xi.\end{eqnarray}
  It follows from Lemma 2.1  that \begin{eqnarray*}\quad \quad \,\;\mbox{Tr} \big((\lambda I - \mathbf{A})^{-1}\big)&=&
    (\lambda -\tau |\xi|^2)^{-1}
  (\lambda-(2\tau+\mu)|\xi|^2)^{-1}  \\
   &&   \times \big[ n\lambda-
  \big((2n-1)\tau +(n-1)\mu\big)|\xi|^2\big],
  \quad \quad \mbox{for}\;\; n=1,2,3,\cdots.\nonumber\end{eqnarray*}
   Thus, for any $W\subset \Omega$.  \begin{eqnarray} \label{3.9}  \quad \quad \,\mbox{Tr}\bigg(e^{-tP}\big|_{W}\bigg)&=&\int_{W}\left\{
    \frac{1}{(2\pi)^n} \int_{{\Bbb R}^n} e^{i\langle x-x, \xi\rangle} \bigg[\frac{1}{2\pi i} \int_{\mathcal{C}} e^{-t\lambda} \big(\mbox{Tr} \big((\lambda I - A)^{-1}\big)\big) d\lambda\bigg] d\xi\right\}dx\\
   &=& \int_{W}\left\{\frac{1}{(2\pi)^n} \int_{{\Bbb R}^n} \bigg[\frac{1}{2\pi i} \int_{\mathcal{C}} e^{-t\lambda}
    \big(\lambda -\tau |\xi|^2\big)^{-1}
  \big(\lambda-(2\tau+\mu)|\xi|^2\big)^{-1} \right. \nonumber\\
   &&  \left. \times \bigg( n\lambda-
  \big((2n-1)\tau +(n-1)\mu\big)|\xi|^2\bigg) d\lambda\bigg]
    d\xi\right\}dx.\nonumber\end{eqnarray}
 Applying the residue theorem \cite{Ahl}, we obtain
 \begin{eqnarray} \label{3.10} && \frac{1}{2\pi i} \int_{\mathcal{C}} e^{-t\lambda}  \frac{\big( n\lambda- \big((2n-1)\tau+(n-1)\mu\big) |\xi|^2\big)}{\big(\lambda -\tau|\xi|^2\big) \big(\lambda -(2\tau+\mu)|\xi|^2\big)} d\lambda\\
   && \quad \quad  = \frac{n \tau |\xi|^2 -\big((2n-1)\tau +(n-1)\mu\big) |\xi|^2}{\tau|\xi|^2 -(2\tau+\mu)|\xi|^2}e^{-t \tau|\xi|^2} \nonumber\\
   &&\quad \quad  + \frac{n (2\tau+\mu)  |\xi|^2 -\big((2n-1)\tau +(n-1)\mu\big) |\xi|^2}{(2\tau+\mu)|\xi|^2 -\tau|\xi|^2}e^{-t (2\tau+\mu)|\xi|^2}\nonumber\\
 && \quad \quad = (n-1) e^{-t\tau|\xi|^2} +  e^{-t(2\tau+\mu)|\xi|^2}.\nonumber\end{eqnarray}
 From (\ref{3.9}) and (\ref{3.10}), we get
  \begin{eqnarray} \label{3.11} && \mbox{Tr}\bigg(e^{-tP}\big|_{W}\big)= \int_{W}\left\{
    \frac{1}{(2\pi)^n} \int_{{\Bbb R}^n} \bigg[(n-1) e^{-t\tau|\xi|^2} +  e^{-t(2\tau+\mu)|\xi|^2}\bigg]
         d\xi\right\}dx\\
         &=& \int_{W} \bigg[ \frac{n-1}{(4\pi \tau t)^{n/2}} +  \frac{1}{(4\pi (2\tau+\mu) t)^{n/2}} \bigg] dx\nonumber\\
         &=& \bigg[  \frac{n-1}{(4\pi \tau t)^{n/2}} +  \frac{1}{(4\pi (2\tau+\mu) t)^{n/2}}\bigg]|W|.
                  \nonumber\end{eqnarray}
    More importantly, (\ref{3.11}) is still valid for any $n$-dimensional normal coordinate patch $W$ covering a patch of $W\cap \partial \Omega$.

 It remains to consider the boundary asymptotics. Let $x=(x'; x_n)$ be local coordinates for $\Omega$ near $\partial \Omega$.
 If $\mathfrak{E}$ is a local frame on $\partial \Omega$; extend $\mathfrak{E}$ to an $n$-dimensional local frame in a neighborhood of $\partial \Omega$ by parallel transport along the geodesic normal rays (see, p.$\,$1101 of \cite{LU}).
We will apply an  ``imagine method'',  which stems from Mckean-Singer in \S5 of \cite{MS}, to deal with the case of the boundary. Let $\mathcal{M}=\Omega \cup (\partial \Omega)\cup \Omega^*$ be the (closed) double of $\Omega$, and $Q$ the double to $\mathcal{M}$ of
 the  operator $P$. Define $Q^{-}$ and $Q^{+}$ to be $Q\big|_{C^\infty(\bar \Omega)}$ subject to $\mathbf{u}=0$ and $\frac{\partial \mathbf{u}}{\partial \nu}=0$ on $\partial \Omega$, respectively. The  $\frac{\partial \mathbf{u}}{\partial t}=Q\mathbf{u}$ still has a nice fundamental solution $\mathbf{K}(t, x,y)$ of class $C^\infty [(0,\infty)\times (\mathcal{M} \setminus\partial \Omega)^2 ]
 \cap C^1 ({\mathcal{M}}^2)$, approximable even on $\partial \Omega$, and the fundamental solution $\mathbf{K}^{\mp}(t,x,y)$
 of $\frac{\partial \mathbf{u}}{\partial t}=Q^{\mp}\mathbf{u}$ can be expresses on $(0,\infty)\times \Omega\times \Omega$ as
 \begin{eqnarray} \label{c4-23} {\mathbf{K}}^{\mp} (t,x,y) =\mathbf{K}(t,x,y)\mp \mathbf{K}(t,x,\overset{\ast} {y}),\end{eqnarray}
 $\overset{*} {y}$ being the double of $y\in \Omega$ (see, p.$\,$53 of \cite{MS}). We pick a self-double patch $W$ of $\mathcal{M}$ covering a patch $W\cap \partial \Omega$ endowed (see the diagram on p.$\,$53 of \cite{MS}) with local coordinates $x$ such that $\epsilon>x_n>0$ in $W\cap \Omega$; $\,x_n=0$ on $W\cap \partial \Omega$;
  $\; x_n (\overset{*}{x})=-x_n(x)$; and the positive $x_n$-direction is perpendicular to $\partial \Omega$. This products the following effect that \begin{eqnarray*} \delta_{jk} (\overset{*}{x})&=&- \delta_{jk} (x) \quad \, \mbox{for}\;\;
  j<k=n \;\;\mbox{or}\;\; k<j=n,\\  &=& \delta_{jk} (x) \;\;\mbox{for}\;\; j,k<n \;\;\mbox{or}\;\; j=k=n,\\
    \delta_{jk}(x)&=& 0 \;\; \mbox{for}\;\; j<k=n \;\;\mbox{or}\;\; k<j=n \;\;\mbox{on}\;\; \partial \Omega.\end{eqnarray*}
     For any $W\subset \Omega$, by the previous technique  we see that (\ref{3.11}) still holds. For any $n$-dimensional normal coordinate patch $W$ covering a patch $W\cap \partial \Omega$, it follows from (\ref{3.11})
     that  \begin{eqnarray} \label{3.12}    &&  \int_{W\cap \Omega}\mbox{Tr}\big( \mathbf{K}(t, x, x)\big)dx= \frac{1}{(4\pi \mu t)^{\frac{n}{2}}} \bigg[\int_{W\cap \Omega}
      \bigg[  \frac{n-1}{(4\pi \tau t)^{n/2}} +  \frac{1}{(4\pi (2\tau+\mu) t)^{n/2}}\bigg] dx\\
       && =\bigg[  \frac{n-1}{(4\pi \tau t)^{n/2}} +  \frac{1}{(4\pi (2\tau+\mu) t)^{n/2}}\bigg]|W\cap \Omega|
   \quad \; \; \mbox{for all }\;\; t\ge 0.\nonumber\end{eqnarray}
     Next, for any small $n$-dimensional normal coordinate patch $W$ covering a patch of $W\cap \partial \Omega$, noting that $|x-\overset{\ast} {x}|=x_n-(-x_n)=2x_n$ we find by the method of pseudodifferential operator that
      \begin{eqnarray} \label{aa-3-1-2}
 && \int_{W\cap \Omega}\mbox{Tr}\big( \mathbf{K}(t, x, \overset{*}{x})\big)dx=  \int_0^\epsilon dx_n \int_{W\cap \partial \Omega}
 \frac{dx'}{(2\pi)^n}\int_{{\Bbb R}^{n-1}} e^{i\langle x-\overset{*}{x}, \xi\rangle}\nonumber \\
  &&\quad  \times \bigg[\frac{1}{2\pi i} \int_{\mathcal{C}} e^{-t\lambda} \big(\mbox{Tr} \big((\lambda I - A)^{-1}\big)\big) d\lambda\bigg] d\xi\nonumber\\
 && =  \int_0^\epsilon dx_n \int_{W\cap \partial \Omega}
 \frac{dx'}{(2\pi)^n}\int_{{\Bbb R}^{n-1}} e^{i\langle x-\overset{*}{x}, \xi\rangle} \bigg[(n-1) e^{-t\tau|\xi|^2} +  e^{-t(2\tau+\mu)|\xi|^2}\bigg] d\xi\nonumber\\
  && = \int_0^\epsilon dx_n \int_{W\cap \partial \Omega}  \frac{dx'}{(2\pi)^n}\int_{-\infty}^\infty e^{2ix_n \xi_n}\bigg\{ \int_{{\Bbb R}^{n-1}} e^{\langle 0, \xi'\rangle }  \bigg[ (n-1) e^{-t\tau(|\xi'|^2+|\xi_n^2)} \nonumber\\
  && \quad +  e^{-t(2\tau+\mu)(|\xi'|^2+x_n^2)} \bigg] d\xi' \bigg\}d\xi_n\nonumber\\
  && = \int_0^\epsilon dx_n \int_{W\cap \partial \Omega}  \frac{dx'}{(2\pi)^n}\int_{-\infty}^\infty e^{2ix_n \xi_n}e^{-t\tau \xi_n^2} \bigg\{ \int_{{\Bbb R}^{n-1}} \bigg[ (n-1) e^{-t\tau \sum_{j=1}^{n-1}\xi_j^2}  \bigg] d\xi' \bigg\}d\xi_n\nonumber\\
  && \quad\,+ \int_0^\epsilon dx_n \int_{W\cap \partial \Omega}  \frac{dx'}{(2\pi)^n}\int_{-\infty}^\infty e^{2ix_n \xi_n}e^{-t(2\tau+\mu) \xi_n^2} \bigg\{ \int_{{\Bbb R}^{n-1}} \bigg[ e^{-t(2\tau+\mu) \sum_{j=1}^{n-1}\xi_j^2}  \bigg] d\xi' \bigg\}d\xi_n\nonumber\\
  && = \int_0^\epsilon dx_n \int_{W\cap \partial \Omega} \frac{n-1}{(4\pi \tau  t)^{n/2}}\, e^{-\frac{(2x_n)^2}{ 4\tau t} }dx' + \int_0^\epsilon dx_n \int_{W\cap \partial \Omega} \frac{1}{(4\pi (2 \tau+\mu)  t)^{n/2}}\, e^{-\frac{(2x_n)^2}{ 4(2\tau+\mu) t} }dx' \nonumber\\
   && = \int_0^\infty dx_n \int_{W\cap \partial \Omega} \frac{n-1}{(4\pi \tau  t)^{n/2}}\, e^{-\frac{(2x_n)^2}{ 4\tau t} }dx'- \int_\epsilon^\infty dx_n \int_{W\cap \partial \Omega} \frac{n-1}{(4\pi \tau  t)^{n/2}}\, e^{-\frac{(2x_n)^2}{ 4\tau t} }dx'\nonumber\\
   && \quad \, + \int_0^\infty dx_n \int_{W\cap \partial \Omega} \frac{1}{(4\pi (2\tau+\mu)  t)^{n/2}}\, e^{-\frac{(2x_n)^2}{ 4(2\tau+\mu) t} }dx' \nonumber \\ && \quad - \int_\epsilon^\infty dx_n \int_{W\cap \partial \Omega} \frac{1}{(4\pi (2\tau +\mu) t)^{n/2}}\, e^{-\frac{(2x_n)^2}{ 4(2\tau+\mu) t} }dx'\nonumber\\
   &&= \frac{n-1}{4} \cdot \frac{|W\cap \partial \Omega|}{(4\pi \tau t)^{(n-1)/2}} + \frac{1}{4} \cdot \frac{|W\cap \partial \Omega|}{(4\pi (2\tau +\mu) t)^{(n-1)/2}}\nonumber\\
   && \quad \, -   \int_{W\cap \partial \Omega} \left\{\int_\epsilon^\infty\bigg[\frac{n-1}{(4\pi \tau  t)^{n/2}}\, e^{-\frac{(2x_n)^2}{ 4\tau t} }+ \frac{1}{(4\pi (2\tau +\mu) t)^{n/2}}\, e^{-\frac{(2x_n)^2}{ 4(2\tau+\mu) t}} \bigg]dx_n\right\}dx', \nonumber\end{eqnarray}
  where $\xi=(\xi', \xi_n)\in {\Bbb R}^n$, and $\epsilon>0$ is some fixed real number.
     It is easy to verify that for any  fixed $\epsilon>0$ and any integer $m\ge 1$,
   \begin{eqnarray*} &&\int_\epsilon^\infty \frac{1}{(4\pi \mu t)^{\frac{n}{2}}} e^{-\frac{(2x_n)^2}{4\tau t}}dx_n = o(t^{m-\frac{n}{2}})\quad \; \, \mbox{as} \;\, t\to 0^+,\\
   && \int_\epsilon^\infty \frac{1}{(4\pi (2\tau+\mu) t)^{\frac{n}{2}}} e^{-\frac{(2x_n)^2}{4(2\tau+\mu) t}}dx_n = o(t^{m-\frac{n}{2}})\quad \; \, \mbox{as} \;\, t\to 0^+
   .\end{eqnarray*}
   Combining these we see that for any integer $m\ge 1$,
     \begin{eqnarray}\label{3..15} && \int_{W\cap \Omega}\mbox{Tr}\big( \mathbf{K}(t, x, \overset{*}{x})\big)dx=  \frac{n-1}{4} \cdot \frac{|W\cap \partial \Omega|}{(4\pi \tau t)^{(n-1)/2}} + \frac{1}{4} \cdot \frac{|W\cap \partial \Omega|}{(4\pi (2\tau +\mu) t)^{(n-1)/2}}\\
  && \quad\,\quad \;\,
   +o(t^{m-\frac{n}{2}})\quad \; \mbox{as} \, \; t\to 0^+.\nonumber\end{eqnarray}
     It follows from (\ref{3.11}), (\ref{c4-23}), (\ref{3.12}) and (\ref{3..15}) we obtain
 \begin{eqnarray} \label{a-4-1-3}  && \int_{W\cap \Omega}\mbox{Tr}\big( {\mathbf{K}}^{\mp}(t, x, x) \big) dx = \int_{W\cap \Omega}\mbox{Tr}\big({\mathbf{K}}(t, x, x) \big) dx \\
&&\quad \, \mp \int_{W\cap \Omega}\mbox{Tr}\big( {\mathbf{K}}(t, x, \overset{*}{x})\big)dx = \bigg[  \frac{n-1}{(4\pi \tau t)^{n/2}}  +  \frac{1}{(4\pi (2\tau+\mu) t)^{n/2}}\bigg]|W\cap\Omega| \nonumber\\
&& \quad \,\mp \frac{1}{4}\bigg[(n-1) \frac{|W\cap \partial \Omega|}{(4\pi \tau t)^{(n-1)/2}}+ \frac{|W\cap \partial \Omega|}{(4\pi (2\tau +\mu) t)^{(n-1)/2}}\bigg]\nonumber
  \\   && \quad\,  +O(t^{m-\frac{n}{2}})\quad \; \mbox{as} \, \; t\to 0^+,\nonumber\end{eqnarray}
and hence  (\ref{1-7}) holds.
   $\square$

\vskip 0.30 true cm

\noindent{\bf Remark 3.1.} \ \  Unlike the Laplacian on Riemannian manifold or other elliptic operators (see, for example, \cite{Min}, \cite{CLN}, \cite{Gil}, \cite{MS} or \cite{Liu1}), in the expansion of the heat trace for the Navier-Lam\'{e} operator, there exist only two spectral invariants (i.e., volume and surface area). This originates from the following two facts: (i) the Lam\'{e} coefficients $\tau$ and $\mu$ are constants (so that the Navier-Lam\'{e} operator is an elliptic operator of  constant coefficients in ${\Bbb R}^n$); (ii) the expansion (\ref{1-7}) is an exact formula except for an additional term $O(t^{\frac{n-1}{2}}e^{-\frac{c}{\sqrt{t}}})$, which is exponential decay as $t\to 0^+$.

\vskip 0.35 true cm

 Now, we  use the heat invariants of the Navier-Lam\'{e} spectrum which have been obtained from Theorem 1.1 to finish the proof of Corollary  1.2.

\vskip 0.35 true cm

 \noindent  {\it Proof of Corollary 1.2.} \  By Theorem 1.1, we know that the first two coefficients  $|\Omega|$ and $|\partial \Omega|$ of the  asymptotic expansion in (\ref{1-7}) are
  Navier-Lam\'{e} spectral invariants, i.e., $|\Omega|= |B_r|$ and $|\partial \Omega|=|\partial B_r|$. Thus $\frac{|\partial \Omega|}{|\Omega|^{(n-1)/n}} =
 \frac{|\partial B_r|}{|B_r|^{(n-1)/n}}$. Note that for any $r>0$, $\frac{|\partial B_r|}{|B_r|^{(n-1)/n}} =
 \frac{|\partial B_1|}{|B_1|^{(n-1)/n}}$.  According to the classical isometric inequality (which states that for any bounded domain $\Omega\subset {\Bbb R}^n$ with smooth boundary, the following inequality holds:
 \begin{eqnarray*}\frac{|\partial \Omega|}{|\Omega|^{(n-1)/n}} \ge
 \frac{|\partial B_1|}{|B_1|^{(n-1)/n}}.\end{eqnarray*}
 Moreover, equality obtains if and only if $\Omega$ is a ball, see \cite{Ch1} or p.$\,$183 of \cite{CLN}), we immediately get $\Omega=B_r$.  $\;\; \square$

\vskip 0.28 true cm

\noindent{\bf Remark 3.2.} \ \    By applying the Tauberian theorem (see, for example, Theorem 15.3 of p.$\,$30 of \cite{Kor} or p.$\,$446 of \cite{Fel}) for the first term on the right side of (\ref{1-7}) (i.e., $\sum_{k=1}^\infty e^{-t\lambda_k^{\mp}}=\int_0^\infty e^{-t\eta} dN(\eta)= \big[ \frac{n-1}{(4\pi \tau t)^{n/2}} +  \frac{1}{(4\pi (2\tau+\mu) t)^{n/2}}\big]|\Omega|+o(t^{n/2})$ as $t\to 0^+$), we can easily obtain the Weyl-type law  for the Navier-Lam\'{e} eigenvalues:
\begin{eqnarray}  N(\eta)&=&\max \{k\big| \lambda_k^{\mp} \le \eta\} = \frac{|\Omega|}{\Gamma (\frac{n}{2} +1)} \bigg[ \frac{n-1}{(4\pi \tau)^{n/2} } \\
 && +\frac{1}{(4\pi (2\tau +\mu))^{n/2}}\bigg]\eta^{\frac{n}{2}} +o(\eta^{\frac{n}{2}}), \,\quad \; \mbox{as}\;\; \eta \to +\infty.\nonumber\end{eqnarray}.

 \vskip 0.58 true cm

\vskip 0.98 true cm

\centerline {\bf  Acknowledgments}

\vskip 0.39 true cm
   This
research was supported by SRF for ROCS, SEM (No. 2004307D01)
   and NNSF of China (11171023/A010801). The author was also supported by Beijing Key Laboratory on MCAACI, Beijing Institute of Technology.

  \vskip 1.68 true cm

\vskip 0.32 true cm

\end{document}